\documentclass[final]{siamltex}
\usepackage[dvips]{graphicx}
\usepackage{amsmath}
\usepackage{latexsym}
\usepackage{amssymb}
\usepackage{color}
\usepackage{algorithm,algorithmic}
\title{Numerical methods for differential linear matrix equations via Krylov subspace methods}
\author{M. Hached \thanks{Laboratoire P.  Painlev\'e UMR 8524, UFR de Math\'{e}matiques, Universit\'e des Sciences et Technologies de Lille, IUT A, Rue de la Recherche, BP 179, 59653 Villeneuve d'Ascq Cedex, France, mustapha.hached@univ-lille.fr}  \and K. Jbilou\thanks{L.M.P.A, Universit\'e du Littoral C\^ote d'Opale, 50 rue F. Buisson BP. 699, F-62228 Calais Cedex, France,  jbilou@lmpa.univ-littoral.fr}} 
\date{ }

\newtheorem{rem}{Remark}
\newcommand{\R}{\hbox{ I\hskip -2pt R}}

\begin {document}
\maketitle

\begin{abstract}
In the present paper, we present some numerical methods for computing approximate solutions to some  large differential linear matrix equations. In the first part of this work, we deal with differential generalized Sylvester matrix equations with full rank right-hand sides using a global Galerkin  and a norm-minimization approaches. In the second part, we consider large differential Lyapunov matrix equations with low rank right-hand sides and use the extended global Arnoldi process to produce low rank approximate solutions. We give some theoretical results and present some numerical experiments.
\end{abstract}

{\small {\bf AMS} {\it subject classification}: 65F10}.\\


\section{Introduction}

In this work, we are interested in the numerical solution of two differential linear matrix equations. First, we consider the linear matrix differential equation with a full right-hand side 
\begin{equation}\label{linmatrix}
\left\{
\begin{array}{l}
\dot{X}(t) = \displaystyle \sum_{i=1}^q A_i\, X(t) \, B_i+C,\\
X(t_0)=X_0,\; \; t \in [t_0, \, T_f],
\end{array}
\right.
\end{equation}
where $A_i \in \mathbb{R}^{n \times n}$,  $B_i \in \mathbb{R}^{p
	\times p}$, $i=1,\ldots,q$, $C\, {\rm and}\, X \in \mathbb{R}^{n
	\times p}$, and we assume that the right hand term $C$ is  full rank. 
Differential Sylvester and Lyapunov   matrix equations are particular cases of (\ref{linmatrix}).\\
The second differential matrix equation that will be considered in this paper, is the well known differential Lyapunov matrix equation with  a low rank right hand side
\begin{equation}\label{lyap1}
\left\{
\begin{array}{l}
\dot X(t)=A\,X(t)+X(t)\,A^T+BB^T;\; (DLE) \\
X(t_0)=X_0,\; \; t \in [t_0, \, T_f],
\end{array}
\right.
\end{equation}

\noindent where the matrix $ A \in \mathbb{R}^{n \times n} $ is assumed to be large, sparse and nonsingular and   $ B \in
\mathbb{R}^{n \times p} $  is  a full rank matrix, with $ p \ll n $.  The initial condition $X_0={\widetilde Z}_0 {\widetilde Z}_0 ^T$ is assumed to be a given symmetric and positive low-rank  matrix. \\
The differential linear matrix equations \eqref{linmatrix} and \eqref{lyap1}  play an important role in many areas   such as control,
filter design theory, model reduction problems, differential
equations and robust control problems \cite{abou03,corless}.\\ Notice that the two linear differential matrix equations above can be reformulated as
\begin{equation}
\label{lin1}
{\dot x}(t) = {\cal M} x(t)+ b,
\end{equation}
where $x(t)=vec(X(t))$, the matrix ${\cal M}$ is given by  ${\cal M} = \displaystyle \sum_{1}^q (B_i^T \otimes A_i)$ for the problem \eqref{linmatrix} and ${\cal M} = I \otimes A +A\otimes I$ for \eqref{lyap1},  while the right hand side $b$ is given by $b=vev(C)$ for \eqref{linmatrix} and $b=vev(BB^T)$ for \eqref{lyap1} respectively, and $vec(C)$ is the long vector obtained by stacking the columns of the matrix $C$.  So we can use classical solvers for computing solutions of \eqref{lin1}. \\
The exact solution of \eqref{lin1} is given by
\begin{equation}
\label{lin2}
x(t) = e^{t {\cal M}} x_0 + \displaystyle \int_{t_0}^t e^{(t-\tau){\cal M}} b \, d \tau,
\end{equation}
which can also be expressed as
\begin{equation}
\label{lin3}
x(t)=t \psi_1(t {\cal M}) (b+{\cal M}x_0)+x_0,
\end{equation}
where 
\begin{equation}
\label{lin4}
\psi_1(z)=\displaystyle \frac{e^z-1}{z}.
\end{equation}
However, in the cases for which the matrix ${\cal M}$ is very large, this approach would not be appropriate. \\

\noindent The rest of the paper is organized as follows. In Section 2, we recall the definitions of the Kronecker and the $\diamond$ products with some of their properties that will be of use in this work. In Section 3, we give a numerical method for solving the problem \eqref{linmatrix} by using projections onto matrix Krylov subspaces, based on a Global-Galerkin orthogonality condition. In Section 4, we will be interested in the numerical solution of the Lyapunov differential matrix equation \eqref{lyap1}. The approximate solutions will be obtained via projection onto matrix Krylov subspaces using the extended global Arnoldi algorithm. The last section is devoted to some numerical experiments.

\section{Preliminaries}
\subsection{Definitions}
We begin by recalling  some notations that  will be used in the sequel. We define the inner product $$\left<Y,Z\right>_F = tr(Y^TZ),$$ where $ tr(Y^TZ)$ denotes the trace of the matrix $Y^TZ$ such that $Y,Z \in \R^{n\times p}$. The associated norm is the Frobenius norm denoted by
$\|Z\|_F=\sqrt{\left<Z,Z\right>_F}$.\\

\noindent The matrix product $A \otimes B = [a_{i,j}B]$ denotes the well known  Kronecker product of the matrices $A$ and $B$ which satisfies the following properties:
\begin{enumerate}
	\item  $(A \otimes B)(C \otimes D) = (AC \otimes BD)$.
	\item  $(A\otimes B)^T = A^T \otimes B^T $.
	\item  $(A\otimes B)^{-1} = A^{-1} \otimes B^{-1} $, if $A$ and $B$ are nonsingular square matrices.
\end{enumerate}
\medskip 

\noindent We also use the matrix  product  $\diamond$ defined in \cite{bouyouli} as follows.

\begin{definition}
	Let ${\cal Z} = [Z_1,...,Z_m]$ and ${\cal W} = [W_1,..., W_l]$ be matrices of dimension $n \times mp$ and $n\times lp$ respectively, where $Z_i$ and $W_j$ $(i = 1,..., m$  $j = 1,...,l)$ are 
	$\R^{n\times p}$. Then the $\R^{m\times l}$ matrix ${\cal Z}^T \diamond {\cal W}$ is defined as:
	$$
	{\cal Z}^T \diamond {\cal W}= [\left< Z_i, W_j \right>]_{1\le i \le m; 1 \le j \le l}
	$$
\end{definition}
A block matrix ${\cal W} = [W_1,..., W_l]$  is said to be F-orthonormal if 
	\begin{equation}
	\left< W_i,W_j \right>= \delta_{i,j} = \left\{\begin{array}{ccccc}
	0 & if & i\neq j \\
	1 & if & i= j
	\end{array}\right. \; i,j=1,...,l.
	\end{equation}
which is equivalent to
$${\cal W}^T \diamond {\cal W}=I_l.$$
\noindent The following proposition gives some properties satisfied by the above product.

\begin{proposition}\label{propd}\cite{bouyouli}
	Let $A,\ B,\ C \in \R^{n\times ps}$, $D \in \R^{n\times n}$, $L \in \R^{p\times p}$ and $α \in \R$. Then we have,
	\begin{enumerate}
		\item $(A + B)^T \diamond C = A^T \diamond C + B^T \diamond C$.
		\item $A^T \diamond (B + C) = A^T \diamond B + A^T \diamond C$.
		\item $(\alpha A)^T \diamond C = \alpha(A^T \diamond C)$.
		\item $(A^T \diamond B)^T = B^T \diamond A$.
		\item $(DA)^T \diamond B = A^T \diamond (D^T B)$.
		\item $A^T \diamond (B(L \otimes I_s)) = (A^T \diamond B)L.$
		\item $\| A^T \diamond B \|_F \leq \| A \|_F \| B \|_F$.
	\end{enumerate}
\end{proposition}
\noindent Let $A$ and $V$ be  $n \times n$ and $n \times p$ matrices, respectively, then the matrix (also called the global) Krylov subspace $\mathbb{K}_m(A,V)$ associated to the pair $(A,V)$ is the subspace of ${\mathbb{R}}^{n \times p}$ generated by $V, AV,\ldots,A^{m-1}V$, i.e.,
$$\mathbb{K}_m(A,V)=span\{V,AV,\ldots,A^{m-1}V\}.$$

In the next proposition, we recall the global QR (gQR) factorisation of an $n \times mp$ matrix $Z$. The algorithm of such a matrix factorisation is given in 
	\cite{bouyouli}. 
\begin{proposition}
	\cite{bouyouli}
	Let $\mathcal{Z} = \left [\begin{array}{llll} Z_{1}, & Z_{2}, &\ldots, &Z_{m} \end{array}\right ]$ be an $n \times mp$ matrix  with $Z_i \in \R^{n \times p}$, $i=1,\ldots,m$. Then,
	the matrix $\mathcal{Z}$ can be factored as 
	$${\mathcal Z} = {\mathcal Q}\,(R \otimes I_p), $$ 
	where ${\mathcal Q}=[Q_1,\ldots,Q_m]$ is an $n \times mp$ F-orthonormal matrix satisfying ${\mathcal Q}^T \diamond {\mathcal Q}=I_m$ and $R$ is  an  upper triangular matrix of dimension  $m \times m$.
\end{proposition}

The following proposition will be useful later.

\begin{proposition}\label{prop4}\cite{jbiloumessaoudi}
	Let  ${\cal V}_m=[V_1,\cdots,V_m]$,  be an $n \times mp$ F-orthonormal matrix.  Let $Z=[z_{i,j}]$ and $G=[g_{i,j}]$
	be  matrices of sizes ${m\times r}$ and ${mp\times q}$ respectively, where $r$ and $q$ are any integers. Then we have
	$$\|{\cal V}_m(Z\otimes I_p)\|_F= \|Z\|_F$$ \noindent
	and
	$$\|{\cal V}_m\,G\|_F\le \|G\|_F. $$
\end{proposition}
\section{Global-based Krylov subspace methods for the problem \eqref{linmatrix}}
\subsection{The Global-Galerkin Krylov subspace  method for linear matrix differential equations}\label{galerkin1}
In this section, we consider the differential linear matrix equation \eqref{linmatrix} and will present an iterative projection method to get numerical approximate solutions. 
Let ${\cal A}$ be the linear matrix operator defined as follows
\begin{eqnarray*}
\mathcal{A}:\;\mathbb{R}^{n \times p} &\longrightarrow&
\;\mathbb{R}^{n \times p}\\
X \; &\longrightarrow& \;\displaystyle \sum_{i=1}^q A_i\, X \,
B_i.
\end{eqnarray*}
Notice that  the transpose of the operator $\mathcal{A}$ with respect to
the inner product $\langle .\,,\, . \rangle_F$ is defined as the application mapping $X \in\mathbb{R}^{n \times p}$ 
$\mathcal{A}^T(X) =\displaystyle \sum_{i=1}^q A_i^T\, X \, B_i^T.$\\
Let $V$ be  any $n \times p$ matrix then we define  the matrix Krylov
subspace associated to the pair $({\mathcal{A}},V)$ and an integer
$m$ defined by $${\cal
K}_{m}({\mathcal{A}},V)=span\{V,{\mathcal{A}}(V),\ldots,{\mathcal{A}}^{m-1}(V)\},$$
 Where  ${\mathcal{A}}^i(V)$ is defined recursively as
${\mathcal{A}}^i(V)={\mathcal{A}}({\mathcal{A}}^{i-1} (V))$.
Notice  that the matrix Krylov subspace  ${\cal
K}_{m}({\mathcal{A}},V)$ is a subspace of $\mathbb{R}^{n \times
p}$, which means that if a matrix $Y$ is in ${\cal
K}_{m}({\mathcal{A}},V)$, then we have $Y=\displaystyle \sum_{i=1}^m \alpha_i {\cal A}^{i-1}(V)$ where $\alpha_i \in \mathbb{R}$, $i=1,\ldots,m$. Next, we remind the 
 modified global Arnoldi algorithm \cite{jms} that allows us to construct an
F-orthonormal basis  $V_1,V_2,\ldots,V_m$\, of the matrix Krylov
subspace ${\cal K}_{m}({\mathcal{A}},V)$, i.e. $$\langle
V_i,V_j\rangle_F=\delta_{i,j},\quad\hbox{for}\quad
i,j=1,\cdots,k,$$ where $\delta_{i,j}$ denotes the classical
Kronecker symbol. The algorithm is described as follows.

\begin{algorithm}[h!]
	\caption{The Modified Global Arnoldi algorithm}\label{alg1}
\begin{itemize}
\item {\tt Set} $V_1=V/\|V\|_F$.
\item  {\tt For} $j=1,\ldots, k$
\begin{enumerate}
	\item ${\tilde
V}={\mathcal{A}}( V_j)$,
\item  {\tt for} $i=1,\ldots, j$.
\begin{enumerate}
\item $h_{i,j}=\langle V_i,{\tilde V}\rangle_F$
\item ${\tilde V}={\tilde
V}-h_{i,j}V_i,$
\item {\tt EndFor}
 \end{enumerate}
\item $h_{j+1,j}=\parallel {\tilde V}
\parallel_F,$
\item  $V_{j+1}={\tilde V}/h_{j+1,j}.$
\end{enumerate}
\item {\tt EndFor}.
\end{itemize}
\end{algorithm}

\noindent 
The matrix  ${\widetilde H}_m$ denotes the $(m+1) \times m$
upper Hessenberg matrix whose nonzero entries $h_{i,j}$ are
defined by Algorithm 1 and $H_m$ is the $m \times m$ matrix
obtained from ${\widetilde H}_m$
 by deleting its last row.  The  $n \times mp$ block matrix ${\cal V}_m=[V_1,V_2,\ldots,V_m]$ is
F-orthonormal which means that the matrices $V_1,\ldots,V_m$ are
orthonormal with respect
 to the scalar product $\langle \,.,.\, \rangle_F$ which is equivalent to 
  \begin{equation}
 \label{v1}
 {\cal V}_{m}^T \diamond {\cal V}_{m} = I_p.
 \end{equation}
For the extended global Arnoldi algorithm, we have the following relations
\begin{equation}\label{Arnoldi1}
[{\mathcal{A}}(V_1),\ldots,{\mathcal{A}}(V_m)]={\cal V}_m(H_m
\otimes I_p) + E_{m+1},
\end{equation}
 where $E_{m+1}=h_{m+1,m}\,[0_{n \times p},\ldots,0_{n
\times p},V_{m+1}]$,
 and
\begin{equation}\label{Arnoldi2}
[{\mathcal{A}}(V_1),\ldots,{\mathcal{A}}(V_m)]={\cal
V}_{m+1}({\widetilde H}_m \otimes I_p).
\end{equation}
Starting from an initial guess $X_0(t) \in \mathbb{R}^{n \times p}$
and the corresponding residual $\mathcal{R}_0(t)={\dot X}_0(t)-{\mathcal{A}}(X_0(t))-C$, at step $m$, we define  the  approximation 
$X_m$ as follows
\begin{equation}\label{approx1}
X_m(t)=X_0(t)+Z_m(t)\;\; {\rm with}\;\; Z_m(t) \in {\cal
K}_m({\mathcal{A}},C)
\end{equation}
and
\begin{equation}\label{ortho1}
\mathcal{R}_m(t)={\dot X}_m(t)-{\mathcal{A}}(X_m(t))-C\,  \perp_F \, {\cal
K}_{m}({\mathcal{A}},C).
\end{equation}
The Galerkin condition \eqref{ortho1} is equivalent to
\begin{equation}
\label{ortho2}
{\cal V}_m^T \diamond \mathcal{R}_m(t) = 0.
\end{equation}
The condition \eqref{approx1} can be written as
\begin{equation}
\label{approx2}
X_m(t) = X_0(t)+\displaystyle \sum_{i=1}^m y_m^{(i)}(t) V_i =X_0(t)+ {\cal V}_m (y_m(t) \otimes I_p),
\end{equation}
where $y_m(t)$ is a vector of $\mathbb{R}^m$ and $y_m^{(i)}(t)$ is the $i$-th component of $y_m(t)$. 
Therefore, the residual $\mathcal{R}_m(t)$ can be expressed as
\begin{eqnarray*}\label{res1}
\mathcal{R}_m(t) &=& {\dot X}_m(t)-{\cal A}(X_m(t))-C\\
& = & {\dot X}_0(t)-{\cal A}(X_0(t))+{\cal V}_m  (\dot{y}_m(t) \otimes I_p)-{\cal A}({\cal V}_m  (y_m(t) \otimes I_p))-C\\
& = & {\cal V}_m  (\dot{y}_m(t) \otimes I_p)- {\cal A}( \displaystyle \sum_{i=1}^p y_m^{(i)}(t) V_i) +{\cal R}_0(t)\\
& = &  {\cal V}_m  (\dot{y}_m(t) \otimes I_p)- ( \displaystyle \sum_{i=1}^p y_m^{(i)}(t) {\cal A}(V_i))+{\cal R}_0(t)\\
& =&  {\cal V}_m  (\dot{y}_m(t) \otimes I_p)- [{\cal A}(V_1), {\cal A}(V_2),\ldots,{\cal A}(V_m)] (y_m(t) \otimes I_p))+{\cal R}_0(t).
\end{eqnarray*}
Using the relation \eqref{Arnoldi1}, it follows that
$$\mathcal{R}_m(t)= {\cal V}_m  (\dot{y}_m(t) \otimes I_p)- {\cal V}_m (H_m \otimes I_p) (y_m(t) \otimes I_p))-E_{m+1}(y_m(t) \otimes I_p)+{\cal R}_0(t).$$
On the other hand, $E_{m+1}$ can be expressed as $E_{m+1} =h_{m+1,m}V_{m+1} [0,0,\ldots,I]$ which can be written as $E_{m+1} =h_{m+1,m}V_{m+1} (\widetilde E \otimes I_p)$. Then we get a new expression of the residual given by
\begin{equation}
\label{res2}
\mathcal{R}_m(t) = {\cal V}_m  (\dot{y}_m(t) \otimes I_p)- {\cal V}_m (H_m y_m(t) \otimes I_p)-h_{m+1,m}V_{m+1} ({\widetilde E}y_m(t) \otimes I_p)+{\cal R}_0(t).
\end{equation}
Using the properties of the $\diamond $ product given in Proposition \ref{propd} and the fact that ${\cal V}_m^T \diamond V_{m+1}=0$, the F-orthogonality condition \eqref{ortho2} reduces to the low dimensional linear differential system of equations
\begin{equation}
\label{sol11}
{\dot y}_m(t)=H_my_m(t)+c_m(t),
\end{equation}
where $c_m =- {\cal V}_m^T \diamond {\cal R}_0(t)\in \mathbb{R}^m$.\\
The solution of the ODE \eqref{sol11} is given by
\begin{equation}
\label{sol12}
y_m(t) = e^{t\, H_m}y_m(0) +\displaystyle \int_0^t e^{(t-\tau)H_m}\, c_m(\tau) d\tau.
\end{equation}
A simple way to compute approximate solutions $x_n \approx x(t_n)$  is to use Euler's method defined as follows  
\begin{equation}
\label{euler}
x_{n+1}=e^{hH_m}x_n+h \psi_1(h H_m)c(t_n),
\end{equation}
where $t_n=nh$ and $h$ is a stepsize and the function $\psi_1$ is defined by \eqref{lin4}.

\noindent In the next algorithm, we summarize the main steps of the global-Galerkin 
\begin{algorithm}
\caption{The global-Galerkin algorithm}\label{alg2}
\begin{enumerate}
\item Choose a tolerance $\varepsilon$ and a maximum number of Arnoldi iteration $m_{max}$.

\item  Compute   $\beta=||{\cal R}_0(t)||_F$,  and $V_1= {\cal R}_0(t)/\beta$.
\item {\tt For} $m=1$ ... $m_{max}$
\begin{enumerate}
\item Construct  the $F$-orthonormal basis $V_1,V_2,\ldots, V_m$ by Algorithm \ref{alg1}.
\item Determine $y_m$ as solution of the problem  the problems \eqref{sol12}.
\item Compute  the residual norm   $||{\cal R}_m(t)||_F$.
 \item If  $||{\cal R}_m(t)||_F<\varepsilon$ stop, else,  Goto (a).
 \end{enumerate}
\item {\tt EndFor}.
\item Compute the approximation $X_m=X_0+{\cal V}_m ( y_m \otimes I_p)$.
\end{enumerate}
\end{algorithm}

\section{Global projection methods  for large differential Lyapunov equations  with low-rank right-hand sides}
In this section, we consider the following large scale differential Lyapunov equation \eqref{lyap1}.
%
Differential Lyapunov  equations play a fundamental role in many topics  such as control,
 model reduction problems, differential equations and robust control problems \cite{abou03,corless}. We notice that, as the problem is large and square, we cannot apply the methods developed in Section 3. \\  The expression of the exact solution is given by 
 \begin{equation}\label{solexacte}
 X(t)=e^{(t-t_0)A}X_0e^{(t-t_0)A^T}+\int_{t_0}^t  e^{(t-\tau)A}\, B B^Te^{(t-\tau)A^T}\, d\tau.
 \end{equation}
 In this section, we consider low-rank approximate solutions to the exact solution $X$ using the global or the extended global Arnoldi process \cite{druskin,heyouni,heyouni1,simoncini}.

 \subsection{Projecting  and using the extended global Arnoldi process}
 We will consider extended global Krylov subspaces associated to the pair $(A,B)$ and defined as follows
 \begin{equation}\label{kry1}
 \mathbb{K}_m(A,B)=span(A^{-m},\ldots,A^{-1}B,B,AB,\ldots,A^{m-1}B).
 \end{equation}
 Notice that
 $$ \mathbb{K}_m(A,B) = {\cal K}_{m}(A,B)\, + \, {\cal K}_m(A^{-1},A^{-1}B), $$
 where ${\cal K}_m(A,B)$ is the global Krylov subspace associated to the pair $(A,B)$. To compute an F-orthonormal basis of $ \mathbb{K}_m(A,B) $, we can use the extended global Arnoldi algorithm defined as follows \cite{heyouni}
 \begin{algorithm}[h!]
 	\caption{The extended global Arnoldi algorithm}\label{extalgo}
 	\begin{enumerate}
 		\item Compute the global  QR decomposition: 
 		$[B, A^{-1}B]= V_1 (R \otimes I_p)$ 
 		\item {\tt For} $j=1,\ldots,m$
 		\begin{enumerate}
 			\item Set $V_j^{(1)}$: the first $p$ columns of $V_j$ and  $V_j^{(2)}$: the  second  $p$ columns of $V_j$,
 			\item Set ${\mathbb V}_j= [{\mathbb V}_{j-1}, V_j ]$ and 
 			$U=[AV_j^{(1)}, A^{-1}V_j^{(2)} ]$,
 			\item F-orthogonalize $U$ w.r. to ${\mathbb V}_j$ to get $V_{j+1}$, i.e.
 				\item {\tt For} $i=1,2,\ldots,j$
 					\begin{enumerate}
 					 \item $H_{i,j} =V_i^T \diamond U$,
 					 \item $U=U-V_i (H_{i,j} \otimes I_p)$
 					 \end{enumerate}
 					\item {\tt EndFor}
 			\end{enumerate}
 		\item Compute the QR decomposition $U=V_{j+1} (H_{j+1,j} \otimes I_p)$
 		\item {\tt EndFor}
 		\end{enumerate}
 \end{algorithm}

 If the upper triangular $2 \times 2$ matrices $H_{j+1,j}$ ($j=1,\ldots,m$) are full rank, then Algorithm \ref{extalgo} computes an F-orthonormal basis of the global extended Krylov subspace $ \mathbb{K}_m(A,B)$, the obtained $n \times 2mp$ matrix $\mathbb{V}_m=[V_1,\ldots,V_m]$ is F-orthonormal
 $$\mathbb{V}_m^T \diamond \mathbb{V}_m =I_{2p}.$$
 Let ${\mathbb T}_m= \mathbb{V}_m^T \diamond (A\mathbb{V}_m)=[T_{i,j]}]$ with $T_{i,j}=V_i^T \diamond (AV_j) \in \mathbb{R}^{2 \times 2}$, $i,j=1,\ldots,m$. Then it can be shown that ${\mathbb T}_m$ is a $2m \times 2m$ upper block Hessenberg matrix whose elements can be obtained from the matrix-coefficients $H_{i,j}$ computed by 
 the extended global Arnoldi algorithm. Let ${\widetilde {\mathbb T}}_m = \mathbb{V}_{m+1}^T \diamond (A \mathbb{V}_m)$, then ${\mathbb T}_m$ can be obtained from  ${\widetilde {\mathbb T}}_m$ by deleting the last 2 rows of ${\widetilde {\mathbb T}}_m$.\\
 
 \noindent We have the following algebraic relations \cite{heyouni}.
 \begin{eqnarray}\label{rel2}
 A\mathbb{V}_m &=& \mathbb{V}_{m+1} ({\widetilde {\mathbb T}}_m \otimes I_p)\\ & = & \mathbb{V}_{m} ({{\mathbb T}}_m
  \otimes I_p) + V_{m+1} T_{m+1,m}(E_m^T \otimes I_p),
 \end{eqnarray}
 where $E_m^T=[0,0,\ldots,I_2]$ the matrix of the last two rows of the identity matrix $I_{2m}$.\\

\noindent Let $X_m(t)$ be the desired low-rank approximate solution  given as 
 \begin{equation}\label{approxL}
 X_m(t) = {\mathbb V}_m (Y_m(t) \otimes I_p)  {\mathbb V}_m^T,\; t\in [t_0,\,T_f], 
 \end{equation}
where $Y_m(t) \in \mathbb{R}^{2m \times 2m}$,  satisfies the Petrov-Galerkin orthogonality condition
 \begin{equation}
 \label{galerkin}
 {\cal V}_m^T R_m(t) {\cal V}_m =0,\; t \in [t_0,\; T_f],
 \end{equation}
 where ${\cal R}_m(t)$ is the residual $ {\cal R}_m(t) = \displaystyle {\dot X}_m(t)-A\,X_m(t)-X_m(t)\,A^T- BB^T $.  Then, from \eqref{approxL} and \eqref{galerkin}, we obtain the low dimensional differential Lyapunov  equation
 \begin{equation}\label{lowlyap}
 \displaystyle {\dot Y}_m(t)- \mathbb{T}_m\,Y_m(t)-Y_m(t)\,\mathbb{T}_m^T  - B_mB_m^T=0, \; t\in [t_0,\,T_f], 
 \end{equation}
 with   $ { B}_m= {\cal V}_m^T \diamond B $. Notice that 
 \begin{equation}
 \label{r}
 [B,A^{-1}B]=V_1(R \otimes I_p),\; {\rm and}\; B_m=r_{1,1}e_1^{(2m)},
  \end{equation}
  with $R=[r_{i,j}]$, $1\le i,j \le 2$ and  $e_1^{(2m)}$ is the first vector of the canonical basis of $ \mathbb{R}^{2m}$.\\The low-dimensional differential Lyapunov equation \eqref{lowlyap} will be solved by using some classical linear  differential equations solvers.\\
 In order to limit the computational effort, we give an upper of the norm of the residual that will allow to stop the iterations without explicitly forming $X_m(t)$ which will be given only at the end of the process.
 \begin{theorem} 
 Let $\mathcal{R}_m(t)$ be the residual obtained at step $m$, then we have
 \begin{equation}\label{normres}
\Vert \mathcal{R}_m(t) \Vert  \le  \sqrt{2}\, \Vert T_{m+1,m} E_m^T Y_m(t) \Vert_F,\; t\in [t_0,\,T_f].
 \end{equation}
 \end{theorem}

\begin{proof}
Using  \eqref{rel2} and \eqref{approxL}, the residual $\mathcal{R}_m(t)=\displaystyle {\dot X}_m(t)-A\,X_m(t)-X_m(t)\,A^T- BB^T  $	is expressed as 
\begin{eqnarray*}
{\mathcal R}_m(t) &= &{\mathbb V}_m ({\dot Y}_m(t) \otimes I_p)  {\mathbb V}_m^T-\left[ \mathbb{V}_{m} ({{\mathbb T}}_m
\otimes I_p) + V_{m+1} (T_{m+1,m}E_m^T \otimes I_p) \right ] (Y_m(t) \otimes I_p)  {\mathbb V}_m^T \\
&- &{\mathbb V}_m (Y_m(t) \otimes I_p)  \left[ \mathbb{V}_{m} ({{\mathbb T}}_m
\otimes I_p) + V_{m+1} (T_{m+1,m}E_m^T \otimes I_p)   \right]^T-BB^T.
\end{eqnarray*}
Therefore, using the fact that $Y_m$ is solution of the low dimensional differential problem \eqref{lowlyap}, the residual can be expressed as follows
\begin{equation*}
\label{res3}
{\mathcal R}_m(t) = \mathbb{V}_{m+1}\,  \left (
\left[ 
\begin{array}{cc}
0& Y_m(t)E_mT_{m+1,m}^T\\
T_{m+1,m} E_m^T Y_m(t) & 0
\end{array}
\right] \otimes I_p  \right)\, {\mathbb V}_{m+1}^T.
\end{equation*}
Therefore, applying Proposition \ref{prop4}, we get for any $t \in [t_0,\; T_f]$ the following upper bound 
\begin{equation}
\label{resm}
\Vert {\mathcal R}_m(t) \Vert_F^2 \le 2\, \Vert T_{m+1,m} E_m^T Y_m(t) \Vert_F^2.
\end{equation}
\end{proof}
 
\subsection{Solving the low dimensional differential Lyapunov equation} 
 We have now  to solve the  low dimensional differential Lyapunov equation \eqref{lowlyap} by some integration method such as the well known  Backward Differentiation  Formula (BDF). We notice that BDF is  especially used for the solution of stiff differential equations.\\  At each time $t_k$, let  $Y_{m,k}$ denote the approximation of $Y_m(t_k)$, where $Y_m$ is a  solution of  (\ref{lowlyap}).  Then, the new approximation $Y_{m,k+1}$ of  $Y_m(t_{k+1})$ obtained at step $k+1$ by ↨$\l$-step BDF is defined  by the implicit relation 
 \begin{equation}
 \label{bdf}
 Y_{m,k+1} = \displaystyle \sum_{i=0}^{l-1} \alpha_i Y_{m,k-i} +h_k \beta {\mathcal F}(Y_{m,k+1}),
 \end{equation} 
 where $h_k=t_{k+1}-t_k$ is the step size, $\alpha_i$ and $\beta_i$ are the coefficients of the BDF method as listed  in Table \ref{tab1} and ${\mathcal F}(X)$ is  given by 
 $${\mathcal F}(Y)= {\mathbb T}_m\,Y+Y\,{\mathbb T}_m^T+\,B_m\,B_m^T.$$
 
 \begin{table}[h!!]
 	\begin{center}
 		\begin{tabular}{c|cccc} 
 			\hline
 			$l$ & $\beta$ &$\alpha_0$ & $\alpha_1$ & $\alpha_2$ \\
 			\hline
 			1 & 1 & 1 & &\\
 			2 & 2/3 & 4/3& -1/3 &\\
 			3 & 6/11 & 18/11 & -9/11 & 2/11\\
 			\hline
 		\end{tabular}
 		\caption{Coefficients of the $l$-step BDF method with $l \le 3$.}\label{tab1}
 	\end{center}
 \end{table}
 \noindent The approximate $Y_{m,k+1}$ solves the following matrix equation
 	\begin{equation*}
 	-Y_{m,k+1} +h_k\beta ({\mathbb T}_m Y_{m,k+1} + Y_{m,k+1} {\mathbb T}_m^T+ B_m B_m^T) + \displaystyle \sum_{i=0}^{l-1} \alpha_i Y_{m,k-i} = 0,
 	\end{equation*}
 which can be written as the following  algebraic Lyapunov matrix equation
 
 \begin{equation}
 \label{lyapbdf}
 {\cal T}_m\, Y_{m,k+1}  + \,Y_{m,k+1} {\cal T}_m^T+ {\cal B}_{m,k}\, {\cal B}_{m,k}^T   =0.
 \end{equation}
 We assume  that at each time $t_k$, the approximation $Y_{m,k}$ is  factorised  as a low rank product  $Y_{m,k}\approx Z_{m,k} {Z_{m,k}}^T$, where $Z_{m,k} \in \mathbb{R}^{n \times m_k}$, with $m_k \ll n$. In that case, the coefficient matrices appearing in \eqref{lyapbdf} are given by 
$${\cal T}_m= h_k\beta {\mathbb T}_m -\displaystyle \frac{1}{2}I \;  \mbox{and} \; {\cal B}_{m,k+1}=[\sqrt{h_k\beta} B_m, \sqrt{\alpha_0}Z_{m,k},\ldots,\sqrt{\alpha_{l-1}} Z_{m,k+1-l}].$$
 The  Lyapunov matrix  equation \eqref{lyapbdf} can be solved by applying direct methods based on Schur decomposition such as the Bartels-Stewart algorithm \cite{bartels,gnv}. \\

In the following algorithm, we summarise the main steps of the extended global Arnoldi method for solving the differential Lyapunov matrix equation \eqref{lyap1}.
\begin{algorithm}[h!]
	\caption{The extended global Arnoldi for differential Lyapunov equations (EgAdl)}
	\begin{enumerate}
		\item Inputs: $A$, $B$ a maximum number of extended Arnoldi iteration $m_max$ and a tolerance $tol$.
		\item Apply the extended global Arnoldi Algorithm to the pair $(A,B)$ to get an F-orthonormal matrix $\mathbb{V}_m=[V_1,\ldots,V_m]$ and the upper block Hessenberg matrix $\mathbb{T}_m$.
		\item Solve the low dimensional problem \eqref{lowlyap} by the BDF method. 
		\item If $\mathcal{R}_m(t) < tol$ stop and  compute the obtained approximate solution.
	\end{enumerate}
\end{algorithm}

 \subsection{Using the approximation of the exponential of a matrix}
In this subsection, we will see how to use the expression \eqref{solexacte} to get low rank approximate solutions.  It is known \cite{higham,horn} that for any square matrix $A$, we have the Cauchy's integral representation
\begin{equation}
\label{cauchy1}
f(A)=\displaystyle \frac{1}{2\pi i}\int_{\Gamma}^{} f(\lambda) (\lambda I -A)^{-1}d \lambda,
\end{equation}
where $f$ is an analytic function on and inside a closed contour $\Gamma  \subset \mathbb{C}$ that encloses the spectrum $\sigma(A)$. A very important topic consists in approximating see \cite{benzi,gallopoulos,hoc,saad2}.
\begin{equation}
\label{cauchy2}
f(A)B= \displaystyle \frac{1}{2\pi i}\int_{\Gamma}^{} f(\lambda) (\lambda I -A)^{-1}Bd \lambda.
\end{equation}
On the other hand, using the global Arnoldi algorithm, we can show \cite{jms} that 
\begin{equation*} 
(\lambda I -A)^{-1}B  \approx   {\cal V}_m \left ( (\lambda I-{\cal H}_m)^{-1} \beta e_1 \otimes I_p)\,  \right),
\end{equation*}
where ${\cal V}_m$ is the F-orthonormal matrix obtained from the global Arnoldi process applied to the pair $(A,V)$ and ${\cal H}_m ={\cal V}_m^T \diamond (A{\cal V}_m)$. Then 
\begin{equation}
\label{cauchy3}  
(\lambda I -A)^{-1}B  \approx \beta {\cal V}_m \left ( (\lambda I-{\cal H}_m)^{-1}  \otimes I_p)\, \tilde E_1 \right),
\end{equation}
where $\tilde E_1= e_1 \otimes I_p$ and $\beta = \Vert B \Vert_F$. Therefore, if the contour $\Gamma$  contains also  the spectrum of ${\cal H}_m$, (which is the case for example  if we choose the countour of field of values of the matrix $A$) we get  
\begin{equation}\label{eqq0}
f(A)B \approx \beta {\cal V}_m \displaystyle \frac{1}{2\pi i}\int_{\Gamma}^{} f(\lambda) \left ( (\lambda I-{\cal H}_m)^{-1}  \otimes I_p)\, \tilde E_1 \right) d \lambda,
\end{equation}
which can be written as
\begin{equation}\label{eqq1}
f(A)B \approx \beta {\cal V}_m \displaystyle \frac{1}{2\pi i}\int_{\Gamma}^{} f(\lambda) \left ( (\lambda I-{\cal H}_m)^{-1}  \otimes I_p)\, \tilde E_1 \right) d \lambda= \beta {\cal V}_m (f({\cal H}_m) \otimes I_p)\, \tilde E_1.
\end{equation}
Using the fact that  $\tilde E_1= e_1 \otimes I_p$, we get 
\begin{equation}\label{eqq2}
f(A)B \approx \beta {\cal V}_m (f({\cal H}_m)e_1 \otimes I_p).
\end{equation}
Notice that using some Kronecker product relations,  the expression \eqref{eqq0} can also be written 
\begin{equation}\label{eqq3}
f(A)B \approx \displaystyle \frac{\beta {\cal V}_m }{2\pi i}\int_{\Gamma}^{} f(\lambda)  \left (\lambda I-\left [{\cal H}_m \otimes I_p \right ]^{-1} \right )  \tilde E_1 d \lambda,
\end{equation}
an then
\begin{equation}
\label{eqq4}  
f(A)B \approx \beta {\cal V}_m f({\cal H}_m \otimes I_p) \tilde E_1.
\end{equation}
The two expressions on the right hand sides in \eqref{eqq2} and \eqref{eqq4}  are the same. 
Applying these results to the function $f(x)=e^x$, we get the approximation to the exponential appearing in the expression of the exact solution \eqref{solexacte}
\begin{equation}
\label{cauchy5}  
e^{(t-\tau)A} B \approx \beta {\cal V}_m (e^{(t-\tau){\cal H}_m} e_1  \otimes I_p).
\end{equation}
Assuming that $X_0=0$, we consider approximations $X_m(t)$ to the solution \eqref{solexacte} as follows
\begin{equation}
\label{approx11} 
X_m(t)=  \int_{t_0}^t Z_m(\tau) {Z_m(\tau)}^T\, d \tau,
\end{equation}
where 
\begin{equation}
\label{approx12}  
Z_m(\tau) = \beta {\cal V}_m (e^{(t-\tau){\cal H}_m} e_1  \otimes I_p).
\end{equation} 
Hence, from \eqref{approx11} and \eqref{approx12}, we get
\begin{equation}
\label{approx3} 
X_m(t)={\cal V}_m (G_m(t) \otimes I_p) {\cal V}_m^T,
\end{equation} 
where 
\begin{equation}
\label{approx4}
G_m(t) = \int_{t_0}^t {\widetilde G}_m(\tau) {\widetilde G}_m(\tau)^T d\tau,
\end{equation} 
with ${\widetilde G}_m(\tau)= \beta e^{(t-\tau){\cal H}_m} e_1$. So, to compute the approximation $X_m(t)$, we need to compute the integral \eqref{approx4} which will be done by using a quadrature formula.\\

\noindent The next theorem states that the matrix function $G_m$ is solution of a low dimensional differential Lyapunov equation.
\begin{theorem}
	The function $G_m$ defined by the relation \eqref{approx4} satisfies the following differential Lyapunov equation,
	\begin{equation}
	\label{low1}
	{\dot G}_m(t)={\cal H}_m G_m(t) + G_m(t) {\cal H}_m^T + \beta^2 e_1 e_1^T.
	\end{equation} 
\end{theorem}
\begin{proof}
	The proof can easily be obtained by deriving   the expression \eqref{approx4}.
\end{proof}
\medskip
\noindent Next, we give a result that allows us the computation of the norm of the residual.
\begin{theorem} \label{tres}
	Let $ X_m(t) = {\cal V}_m(G_m(t) \otimes I_p){\cal V}_m^T $ be the approximation obtained at step $ m $ by the  global 
	Arnoldi  method. Then the residual $ {\cal R}_m(t) $ satisfies
	\begin{equation}
	\label{result2}
	\parallel {\cal R}_m(t) \parallel_F \le  \vert h_{m+1,m} \vert \Vert {\bar G}_m(t) \parallel_2,
	\end{equation}
	where $ {\bar G}_m(t) $ is the last row of $ G_m(t) $.
\end{theorem}
\medskip
\begin{theorem}
	Let $X_m(t)$ be the approximate solution given by \eqref{approx3}. Then we have 
	\begin{equation}
	\label{pertu}
	\displaystyle {\dot X}_m(t)=A\,X_m(t)+X_m(t)\,A^T+{\cal L}_m.
	\end{equation}
	where    ${\cal L}_m =BB^T- L_m-L_m^T$ with $L_m(t)=h_{m+1,m} E_m^T (G_m(t) \otimes I) {\cal V}_m^T$.\\ The error $\mathcal{E}_m(t)=X(t)-X_m(t)$  satisfies the following equation
	\begin{equation}\label{pertu2}
	\displaystyle {\dot {\mathcal E}}_m(t)=A\ {\mathcal E}_m(t)+{\mathcal E}_m(t)A^T-\mathcal{R}_m(t),
	\end{equation}
	and then 
	\begin{equation}
	\label{error3}
	\mathcal{E}_m(t)=e^{(t-t_0)A}\mathcal{E}_{m,0}e^{(t-t_0)A^T}+\int_{t_0}^t e^{(t-\tau)A}\mathcal{R}_m(\tau)e^{(t-\tau)A^T}d\tau,\; t \in [t_0,\, T_f].
	\end{equation}
	where  $\mathcal{E}_{m,0}=\mathcal{E}_m(t_0)$.
\end{theorem}
\medskip
\begin{proof}
	The proof of \eqref{pertu2} is obtained by using the expression \eqref{approx3} of the approximate solution $X_m(t)$ and the relation \eqref{low1}. The expression \eqref{error3} of the error is easily derived by extracting the initial problem \eqref{lyap1} from the expression of the residual ${\cal R}_m(t)$.
\end{proof}

\medskip
\begin{theorem}
	\label{Theoerr2}
	Assume that $X(t_0)=X_m(t_0)$, then we have the following upper bound
	\begin{equation}
	\label{upperbound}
	\parallel \mathcal{E}_m(t) \parallel  \le   \displaystyle  \vert  h_{m+1,m} \vert \, \parallel \bar G_m \parallel_{\infty}  \frac{e^{2(t-t_0)\mu_2(A)}-1}{2 \mu_2(A)},\; \forall t \in [t_0,\, T_f],\\
	\end{equation}
	where $\mu_2(A)=\displaystyle \frac{1}{2} \lambda_{max}(A+A^T)$  is the 2-logarithmic norm   and $\parallel \bar G_m \parallel_{\infty}  =\displaystyle \max_{\tau \in [t_0,\, t]} \parallel \bar G_m(\tau) \parallel$ where 	where $ {\bar G}_m(t) $ is the last row of $ G_m(t) $.
	\end{theorem}
\medskip
\begin{proof}
Using the expression \eqref{error3} of $\mathcal{E}_m(t)$ and the fact that $\parallel e^{tA} \parallel \le e^{\mu_2(A)t}$, we  get 
\begin{eqnarray*}
	\parallel \mathcal{E}_m(t) \parallel &  \le  &  \parallel  h_{m+1,m}  \bar G_m \parallel_{\infty}  \displaystyle \int_{t_0}^t e^{2(t-\tau) \mu_2(A)} d\tau\\
	& \le & \vert h_{m+1,m} \vert \, \parallel \bar G_m \parallel_{\infty} e^{2t\mu_2(A)} \displaystyle \int_{t_0}^t e^{-2\tau \mu_2(A)} d\tau\\
	& =  & \displaystyle  \vert  h_{m+1,m} \vert\, \parallel \bar G_m \parallel_{\infty} \frac{e^{2(t-t_0)\mu_2(A)}-1}{2 \mu_2(A)},
\end{eqnarray*}
which gives the desired result.	
	\end{proof}
	
\medskip
\begin{rem}
Notice that instead of using the global Arnoldi, we can also use the extended global Arnoldi to obtain approximation to $f(A)B$. In this case we have
\begin{equation}
\label{expext1}
e^{(t-\tau)A}B \approx (e^{(t-\tau) {\mathbb T}_m} \otimes I_p){\widetilde B}_m,
\end{equation}
where ${\widetilde B}_m=r_{1,1}e^{(2m)} \otimes I_p$ given by \eqref{r}. Then
\begin{equation}
\label{expext2}
e^{(t-\tau)A}B \approx r_{1,1} (e^{(t-\tau) {\mathbb T}_m} e^{(2m)} \otimes I_p).
\end{equation}

\noindent Therefore, all the relations stated for the global Arnoldi are still valid for the extended block Arnoldi with some variations. From the numerical point of view, the extended global Arnoldi  methodis faster than  global Arnoldi.
	\end{rem}


\section{Conclusion}
We presented in the present paper  different  new approaches for computing approximate solutions to large scale differential differential  matrix equations. These approaches are based on projection onto matrix Krylov subspaces using the globlal and the extended global Arnoldi algorithms. For problems with full rank right hand sides, the problem reduces to the computation of solutions of differential linear systems of equations  by classical methods. In the second part of this work, we considered a differential Lyapunov matrix equation with a decomposed low rank hand sides. The initial problem was projected onto matrix Krylov subspaces to get low dimensional differential Lyapunov equation that is solved by the classical BDF methods. 
 Numerical experiments will be provided to show that both methods are promising for large-scale problems

\end{document}